\documentclass[12pt]{article}
\usepackage{latexsym,amsfonts,amssymb,graphicx,graphics}
\newtheorem{theorem}{Theorem}[section]

\newtheorem{example}[theorem]{Example}

\newtheorem{remark}[theorem]{Remark}
\begin{document}
\title{\Large \bf A relationship between two graphical models of the Kauffman
polynomial}

\author{\small Xian'an Jin\footnote{ \textit{E-mail address:}
xajin@xmu.edu.cn}\\
{\small School of Mathematical Sciences, Xiamen University, Xiamen 361005, P. R. China}\\
} \maketitle \vskip1cm
\begin{abstract}
There are two oriented 4-valent graphical models for the Kauffman
polynomial: one ($HJ$) is obtained by combining Jaeger's formula and
Kauffman-Vogel model for the Homflypt polynomial; the other ($WF$)
is obtained by combining Kauffman-Vogel model for the Kauffman
polynomial and Wu's formula. The main goal of this paper is to
explore the relationship between the two models. We find that there
is an one-to-many correspondence between the terms of $HJ$ model and
the terms of $WF$ model. In addition, we investigate the relation
between trivalent graphical models and 4-valent graphical models of
both the Homflypt and Kauffman polynomials, and observe that there
is a bijection between the terms of the two models.
\end{abstract}



\section{Introduction}

In \cite{Kau1,KV}, Kauffman and Vogel generalized the Homflypt and
Kauffman polynomials from links to 4-valent rigid vertex spatial
graphs. Conversely, an unoriented (resp. oriented) 4-valent plane
graph expansion for the Kauffman (resp. Homflypt) polynomial of
unoriented (resp. oriented) links was obtained, which is implicit in
\cite{KV}. In 1989, Jaeger announced a relation \cite{K}, we shall
call it Jaeger's formula, between the Kauffman polynomial of an
unoriented link diagram and the Homflypt polynomials of some
oriented link diagrams constructed from the unoriented link diagram.
Recently, Wu generalized Jaeger's formula from link diagrams to
4-valent rigid vertex spatial graph diagrams \cite{Wu}. We shall
call it Wu's formula. Note that 4-valent rigid vertex spatial graph
diagrams include 4-valent plane graphs as a special case. In this
paper we shall confine ourselves in 4-valent plane graphs.

\begin{figure}[htbp]
\centering
\includegraphics[width=8cm]{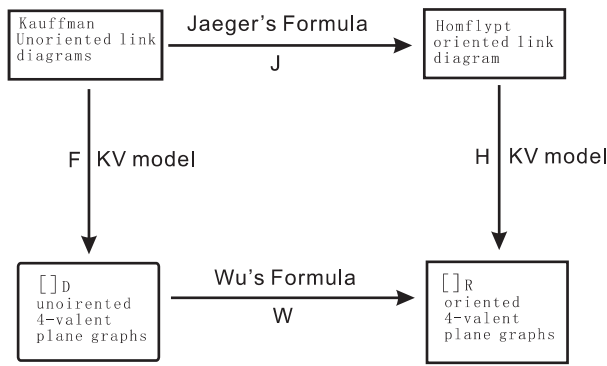}
\renewcommand{\figurename}{Fig.}
\caption{The relation diagram.}\label{cd}
\end{figure}

We illustrate above descriptions by a relation diagram as shown in
Fig. \ref{cd}. According to the diagram, two oriented 4-valent
graphical models for the Kauffman polynomial will be produced: one
($HJ$) is obtained by combining Jaeger's formula and Kauffman-Vogel
model for the Homflypt polynomial; the other ($WF$) is obtained by
combining Kauffman-Vogel model for the Kauffman polynomial and Wu's
formula. The first and main goal of this paper is to explore the
relationship between the two models. We modify the Jaeger's formula
and then find that there is an one-to-many correspondence between
the terms of $HJ$ model and the terms of $WF$ model.

As a result, we actually verified the consistency of the two
Kauffman-Vogel models, Jaeger's formula and Wu's formula. We pointed
that in \cite{Hug}, Huggett found and verified a relationship (i.e.
the replacement of a crossing by a clasp) between a famous
Thistlethwaite's result \cite{MBT} which expresses the Jones
polynomial as a special parametrization of the Tutte polynomial
\cite{Tutte} and a Jaeger's result relating the Homflypt polynomial
with the Tutte polynomial.

Then we investigate the relation between trivalent graphical models
\cite{MOY,Fre,CT} and 4-valent graphical models of both the Homflypt
and Kauffman polynomials. We observed there is a bijection between
the terms of the two models via contracting ``thick" edges of
trivalent plane graphs to obtain 4-valent plane graphs.

\section{Homflypt and Kauffman polynomials}

To proceed rigorously, it is necessary to recall the definition of
the Homflypt and Kauffman polynomials. The Homflypt polynomial was
introduced in \cite{Hom1} and \cite{Hom2}, independently, which is a
writhe-normalization of its regular isotopy counterpart: the $R$
polynomial. Let $L$ be an oriented link diagram. We denote by
$R_L(z,a)\in \mathbf{Z}[z^\pm,a^\pm]$ the $R$ polynomial of $L$.

\begin{center}
{\bf Axioms for the $R$ polynomial}
\end{center}

\begin{figure}[htbp]
\centering
\includegraphics[width=6cm]{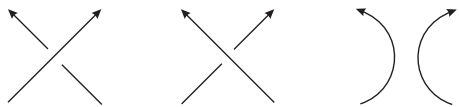}
\renewcommand{\figurename}{Fig.}
\caption{{\footnotesize A skein triple: $L_+$, $L_-$ and
$L_0$.}}\label{skH}
\end{figure}

\begin{enumerate}
\item[(1)] $R_{\circlearrowright}=1$.
\item[(2)] $R_L$ is invariant under Reidemeister moves II and III.
\item[(3)] (the kink formulae) The effect of  Reidemeister move I on $R$ is to
multiply by $a$ or $a^{-1}$ according to the type of Reidemeister
move I:
\begin{eqnarray}
R_{L(+)}=aR_{L},\ \ R_{L(-)}=a^{-1}R_{L}, \label{HR1}
\end{eqnarray}
where $L(+)$ (resp. $L(-)$ denotes diagrams with a positive (resp.
negative) curl and $L$ denotes the result of removing this curl by
Reidemeister move I.
\item[(4)] (the skein relation)
\begin{eqnarray}
R_{L_+}-R_{L_-}=zR_{L_0},
\end{eqnarray}
where $L_+, L_-$ and $L_0$ are link diagrams which are identical
except near one crossing where they are as in Fig. \ref{skH} and are
called a skein triple.
\end{enumerate}

The Alexander-Conway \cite{Alexander,Conway} and Jones \cite{Jones}
polynomials are both special cases of the Homflypt polynomial.

\begin{figure}[htbp]
\centering
\includegraphics[width=8cm]{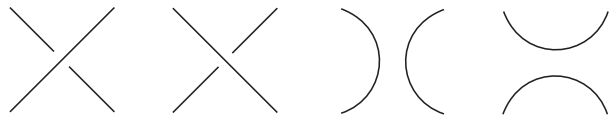}
\renewcommand{\figurename}{Fig.}
\caption{{\footnotesize The quadruple: $L_+$, $L_-$, $L_0$ and
$L_\infty$ (from left to right).}}\label{skK}
\end{figure}

The Kauffman polynomial was introduced in \cite{Kau2-variable}. We
work with its ``Dubrovnik" version \cite{Lik}. Let $L$ be an
unoriented link diagram. We denote by $D_L(z,a)\in
\mathbf{Z}[z^\pm,a^\pm]$ the Dubrovnik (briefly, $D$) polynomial of
$L$. The Dubrovnik polynomial satisfies the following axioms:

\begin{center}
{\bf Axioms for the $D$ polynomial}
\end{center}
\begin{enumerate}
\item[(1)] $D_{\bigcirc}=1$.
\item[(2)] $D_L$ is invariant under Reidemeister moves II and III.
\item[(3)] (the kink formulae) The effect of  Reidemeister move I on $D$ is to
multiply by $a$ or $a^{-1}$ according to the type of Reidemeister
move I:
\begin{eqnarray}
D_{L(+)}=aD_{L},\ \ D_{L(-)}=a^{-1}D_{L}, \label{DR1}
\end{eqnarray}
where $L(+)$ (resp. $L(-)$ denotes diagrams with a positive (resp.
negative) curl and $L$ denotes the result of removing this curl by
Reidemeister move I.
\item[(4)] (the switching formula)
\begin{eqnarray}
D_{L_+}-D_{L_-}=z(D_{L_0}-D_{L_\infty}),\label{Dske}
\end{eqnarray}
where $L_+, L_-, L_0$ and $L_\infty$ are link diagrams which are
identical except near one crossing where they are as shown in Fig.
\ref{skK}.
\end{enumerate}

The Kauffman polynomial is a writhe-normalization of the Dubrovnik
polynomial, which is the generalization of both the Jones polynomial
\cite{Jones} and the BLM-Ho's $Q$ polynomial \cite{BLM,Ho}. The $D$
polynomial specializes to the Kauffman bracket polynomial \cite{KB}
by putting $a=-A^3$ and $z=A-A^{-1}$ \cite{K}.

We point out that, when we mention the Homflypt and Kauffman
polynomials we sometime mean the $R$ and $D$ polynomials,
respectively.

\section{Kauffman-Vogel models}

In this section, we explain the 4-valent graphical models of the $R$
and $D$ polynomials. A graph is \emph{planar} if it can be embedded
in the plane, that is, it can be drawn on the plane so that no two
edges intersect. The embedding of a planar graph is called a
\emph{plane} graph. It is well known that any graph can be embedded
in the 3-dimensional Euclidean space \cite{Bondy}, and such an
embedding is called a \emph{spatial} graph.

A \emph{4-valent graph} is a graph whose each vertex is of degree 4.
We always consider simple closed curves called \emph{free loops} as
special cases of 4-valent graphs, in other words, a free loop is a
graph having one edge and having no vertices.

In \cite{Kau1}, Kauffman defined 4-valent graphs with \emph{rigid
vertices} and introduced the notion of \emph{rigid vertex ambient
isotopy} for 4-valent rigid vertex spatial graphs. In \cite{KV},
Kauffman and Vogel introduced two 3-variable ($A,B,a$) polynomials
for 4-valent rigid vertex spatial graphs in terms of the $R$ and $D$
polynomials, respectively. When $G$ has no vertices, i.e. $G$ is a
link, The two 3-variable polynomials of 4-valent rigid vertex
spatial graph $G$ will specialize to $R$ and $D$ polynomials,
respectively.

Conversely, the $R$ and $D$ polynomials of link diagrams can be
expressed as the sum of such 3-variable polynomials of 4-valent
plane graphs constructed from link diagrams. In \cite{Carp},
Carpentier proved that such 3-variable polynomials can be computed
recursively completely within the category of planar graphs without
resorting to links. Thus we obtain 4-valent plane graphical models
for both the $R$ and $D$ polynomials. Now we give a detailed account
of the two models.

\subsection{$R$ polynomial}

\begin{figure}[htbp]
\centering
\includegraphics[width=1.8cm]{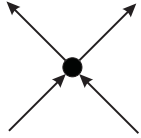}
\renewcommand{\figurename}{Fig.}
\caption{Crossing-like oriented vertex.}\label{crosslike}
\end{figure}

By an \emph{oriented} 4-valent plane graph, we mean a 4-valent plane
graph together with an edge orientation of the graph such that at
each vertex, the four (not necessarily distinct) edges incident with
the vertex are oriented like a crossing of an oriented link diagram
as shown in Fig. \ref{crosslike}.
The 3-variable polynomial $[G]_R=[G]_R(A,B,a)$ for an oriented
4-valent plane graph $G$ can be defined via the following graphical
calculus \cite{KV}.

\begin{center}
{\bf Graphical calculus for the $[]_R$ polynomial}
\end{center}

\begin{enumerate}
\item[(1)]
$[\circlearrowright]_R=1$, where $\circlearrowright$ is a free loop
and its orientation is actually irrelevant.
\item[(2)]
$[G\sqcup \circlearrowright]_R=\delta[G]_R$, where $G\sqcup
\circlearrowright$ is the disjoint union of an oriented 4-valent
plane graph $G$ and $\circlearrowright$, and $\delta={a-a^{-1}\over
A-B}$.
\item[(3)] Let
\begin{eqnarray*}
\lambda&=&{Aa^{-1}-Ba\over A-B}, \\
\theta&=&{B^2a-A^2a^{-1}\over A-B},\\
\eta&=&{B^{3}a-A^3a^{-1}\over A-B}.
\end{eqnarray*}

Then identities as shown in Fig. \ref{id3} hold:
\begin{figure}[htbp]
\centering
\includegraphics[width=10cm]{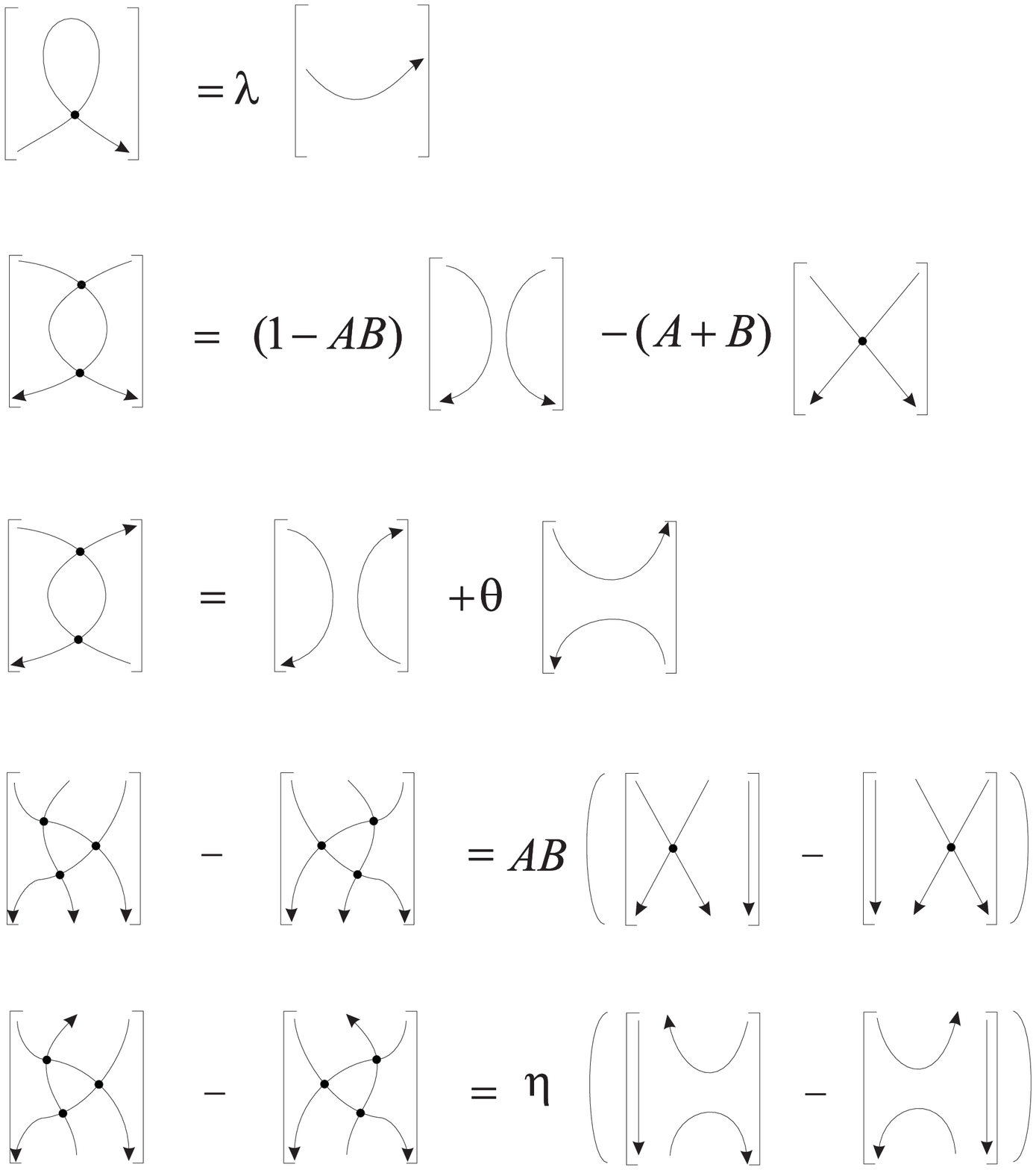}
\renewcommand{\figurename}{Fig.}
\caption{Identities for $[]_R$.}\label{id3}
\end{figure}
\end{enumerate}

The following theorem is implicit in \cite{KV}.

\begin{figure}[htbp]
\centering
\includegraphics[width=4.5cm]{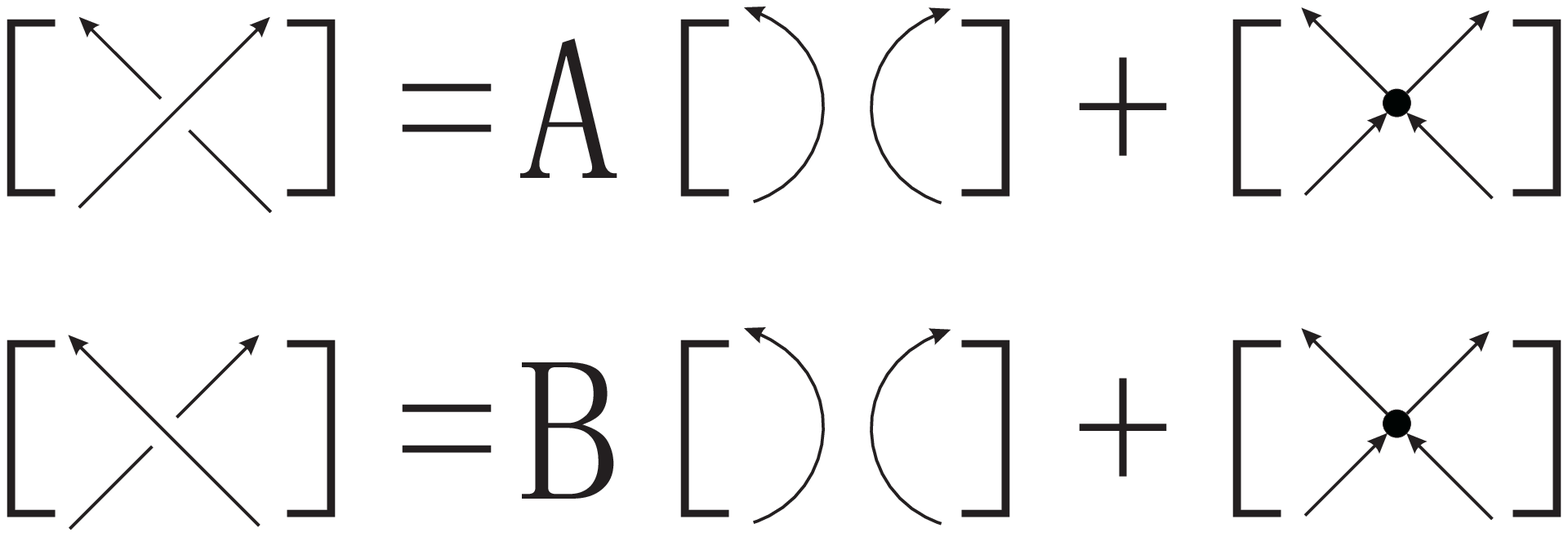}
\renewcommand{\figurename}{Fig.}
\caption{Recursive equation for $R_L$.}\label{He}
\end{figure}

\begin{theorem}\label{yy} Let $L$ be an oriented link diagram. Then
\begin{eqnarray}
R_L(A-B,a)=\sum_{G}A^{i(G)}B^{j(G)}[G]_R,
\end{eqnarray}
where the summation is over all oriented 4-valent plane graphs:
$G$'s, obtained from $L$ by applying to each crossing (positive or
negative) one of the two replacements as shown in Fig. \ref{Hrep},
and $i(G)$ and $j(G)$ are the numbers of positive and negative
crossings of $L$ smoothed to form $G$, respectively.
\begin{figure}[htbp]
\centering
\includegraphics[width=8cm]{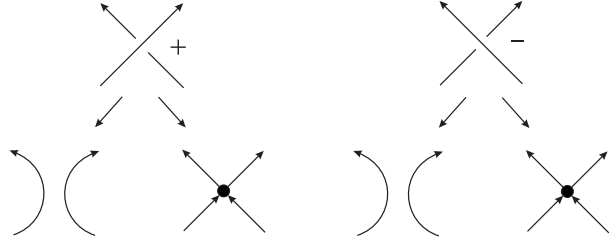}
\renewcommand{\figurename}{Fig.}
\caption{Two types of replacements of an oriented
crossing.}\label{Hrep}
\end{figure}
\end{theorem}
By a little abuse of notations, if we write $[L]=R_L$, then we have
the following recursive equation as shown in Fig. \ref{He}. Note
that the subscript $R$ of $[]_R$ is omitted.

\subsection{$D$ polynomial}

The 3-variable polynomial $[G]_D=[G]_D(A,B,a)$ for an unoriented
4-valent plane graph $G$ can be defined via the following graphical
calculus \cite{KV,Carp}.

\begin{center}
{\bf Graphical calculus for the $[]_D$ polynomial}
\end{center}

\begin{enumerate}
\item[(1)]
$[\bigcirc]_D=1$, where $\bigcirc$ is a free loop.
\item[(2)]
$[G\sqcup \bigcirc]_D=\mu[G]_D$, where $G\sqcup \bigcirc$ is the
disjoint union of an unoriented 4-plane graph $G$ and $\bigcirc$,
and $\mu={a-a^{-1}\over A-B}+1$.
\item[(3)] Let
\begin{eqnarray*}
o&=&{Aa^{-1}-Ba\over A-B}-(A+B), \\
\gamma&=&{B^2a-A^2a^{-1}\over A-B}+AB,\\
\xi&=&{B^{3}a-A^3a^{-1}\over A-B}.
\end{eqnarray*}
Then identities as shown in Fig. \ref{D} hold:
\begin{figure}[htbp]
\centering
\includegraphics[width=10cm]{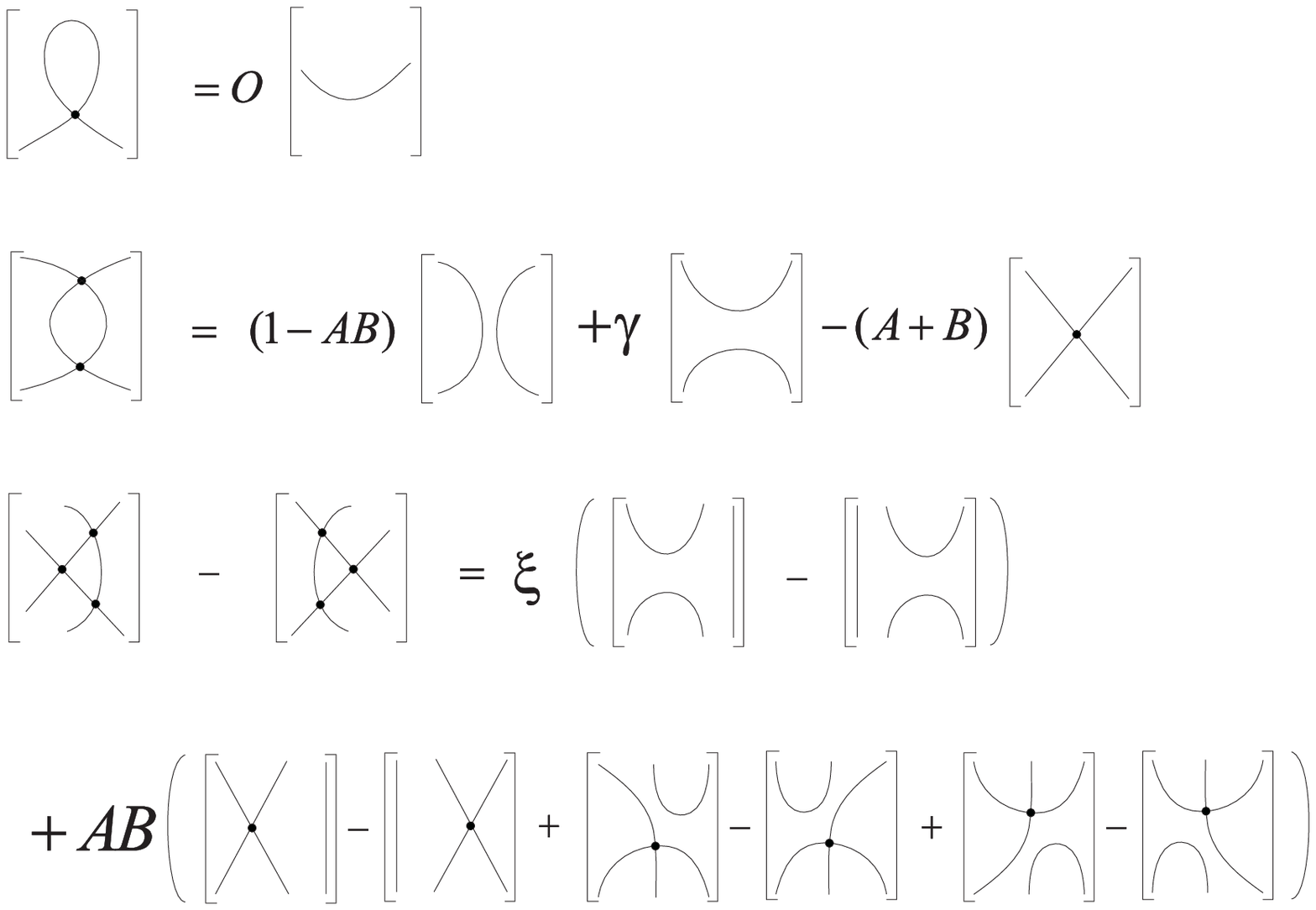}
\renewcommand{\figurename}{Fig.}
\caption{Identities for $[]_D$.}\label{D}
\end{figure}
\end{enumerate}

The following theorem is implicit in \cite{KV}, and explicit in
\cite{Carp}.

\begin{figure}[htbp]
\centering
\includegraphics[width=6cm]{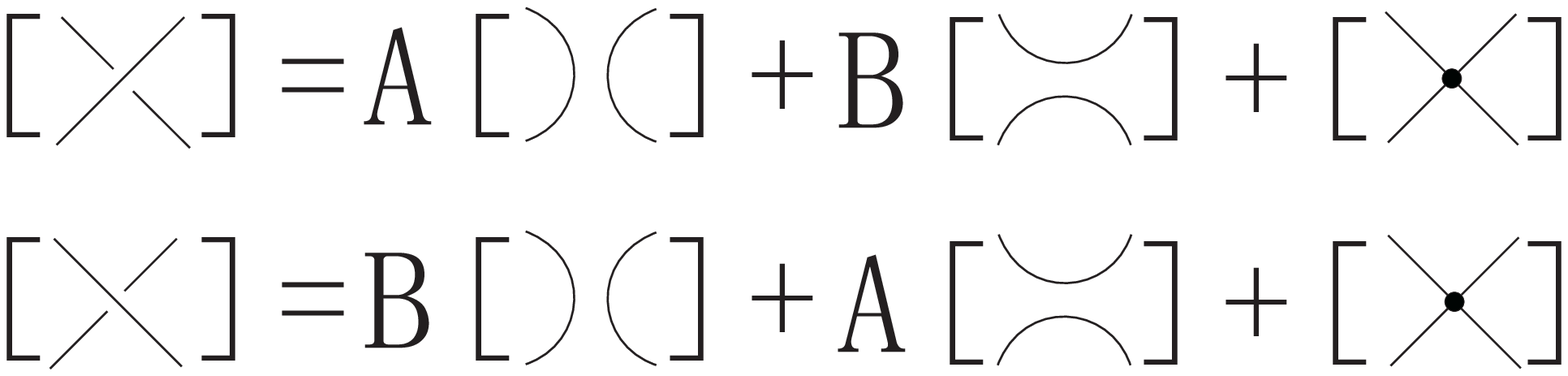}
\renewcommand{\figurename}{Fig.}
\caption{Recursive equation for $D_L$.}\label{Ke}
\end{figure}

\begin{theorem} Let $L$ be an unoriented link diagram. Then
\begin{eqnarray}
D_L(A-B,a)=\sum_{G}A^{i(G)}B^{j(G)}[G]_D,
\end{eqnarray}
where the summation is over all (unoriented) 4-valent plane graphs:
$G$'s, obtained from $L$ by applying to each crossing one of the
three types of replacements as shown in Fig. \ref{Krep}, and $i(G)$
and $j(G)$ are the numbers of crossings of $L$ of $A$-smoothings and
$B$-smoothings used to form $G$, respectively.
\begin{figure}[htbp]
\centering
\includegraphics[width=6cm]{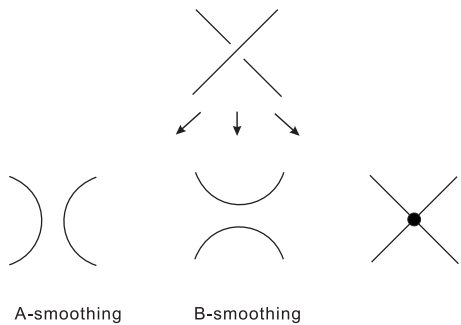}
\renewcommand{\figurename}{Fig.}
\caption{Three types of replacements of an unoriented
crossing.}\label{Krep}
\end{figure}
\end{theorem}
Similarly, if we write $[L]=D_L$, then we have the following
recursive equation as shown in Fig. \ref{Ke}.

\section{Jaeger and Wu's formulae}

Let $L$ be an oriented link diagram, the \emph{rotation number}
(also called Whitney degree, see \cite{KN}, p. 170) rot($L$) of $L$
is equal to the sum of signs for all Seifert circles of $L$ with the
convention that the sign is $+1$ if the circle is counterclockwise
oriented and the sign is $-1$ if the circle is clockwise oriented.
It actually measures the total turn of the unit tangent vector to
the underlying plane curves of the link diagram. The rotation number
of an oriented 4-valent plane graph is thus equal to that of any
oriented link diagram obtained from the graph by converting each
vertex into a crossing.

\subsection{Jaeger's formula}

Jaeger (see \cite{K}, pp. 219-222) established a relation between
Kauffman polynomial and Homflypt polynomial, which expresses the $D$
polynomial of a link diagram as the certain weighted sum of Homflypt
polynomials of some oriented link diagrams obtained by firstly
``splicing" some crossings of the link diagram and then assigning an
orientation. Jaeger's formula can also be found in \cite{Fer,Wu}.
Here we give it a slightly different formulation. Note that in this
paper we use the normalized versions of the $R$ and $D$ polynomials,
while in \cite{K,Fer,Wu}, the authors all dealt with unnormalized
versions.

\begin{figure}[htbp]
\centering
\includegraphics[width=8cm]{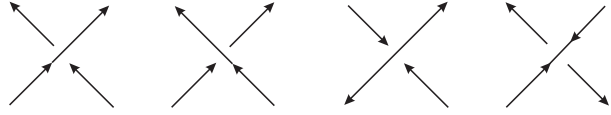}
\renewcommand{\figurename}{Fig.}
\caption{Balanced edge orientation near a crossing: the first is a
positive crossing-like oriented crossing; the second is a negative
crossing-like oriented crossing; the third and the fourth are both
alternatingly oriented called \emph{top outward} crossing and
\emph{top inward} crossing, respectively in \cite{Wu}.}\label{co}
\end{figure}

Let $L$ be an unoriented link diagram. We call a segment of the
diagram between two adjacent crossings an \emph{edge} of $L$. An
edge orientation of $L$ is \emph{balanced} if, at each crossing,
among four (not necessarily distinct) edges around the crossing, two
edges are ``in" and two edges are ``out". Up to rotation, there are
four possible balanced edge orientations near a crossing: two are
crossing-like oriented and the other two are alternatingly oriented.
(See Fig. \ref{co}.)
Denote by $\mathcal{O}(L)$ the set of all balanced edge orientations
of $L$.

\begin{figure}[htbp]
\centering
\includegraphics[width=10cm]{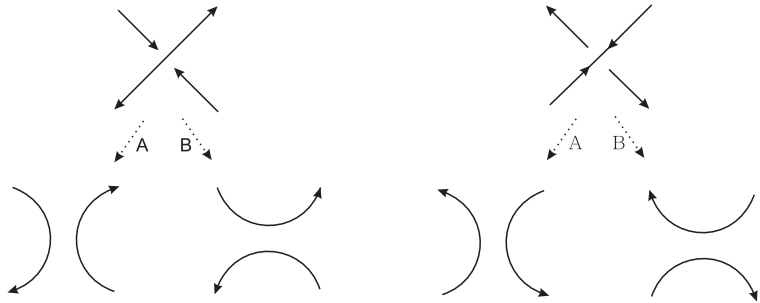}
\renewcommand{\figurename}{Fig.}
\caption{Two ways of smoothing a top outward crossing and a top
inward crossing.}\label{spice-c}
\end{figure}

Given an $o\in \mathcal{O}(L)$, equipping $L$ with $o$, we obtain
$L_o$. To obtain oriented link diagrams from $L_o$ we need to
``splicing" all alternatingly oriented crossings. There are two ways
of smoothing a top outward crossing and a top inward crossing of
$L_o$ as shown in Fig. \ref{spice-c}. A \emph{resolution} $r$ of
$L_o$ is a choice of $A$ or $B$ smoothing for every top outward and
top inward crossing of $L_o$. Denote by $\sum(L_o)$ the set of all
resolutions of $L_o$. Given a $r\in \sum(L_o)$, for each top outward
and top inward crossing $c$ of $L_o$, we define the weight to be

\begin{displaymath}
[L_o,r;c] = \left\{
\begin{array}{ll}
q-q^{-1} & \textrm{if $c$ is a top outward crossing and $r$} \\
         & \textrm{applies $A$ smoothing to $c$,}\\
q^{-1}-q & \textrm{if $c$ is a top outward crossing and $r$}\\
         & \textrm{applies $B$ smoothing to $c$,}\\
0 & \textrm{if $c$ is a top inward crossing.}
\end{array} \right.
\end{displaymath}

We can take the weight of the unspliced crossing $L_o$ to be 1. The
total weight $[L_o,r]$ of the resolution $r$ (applied to $L_o$) is
defined to be the product of weights of all crossings of $L_o$.
Denote by $L_{o,r}$ the oriented link diagram obtained by applying
$r$ to $L_o$.

\begin{theorem}\label{Ja}(\textbf{Jaeger's Formula})  Let $L$ be an unoriented link diagram. Then
\begin{eqnarray*}
D_L(q-q^{-1},a^2q^{-1})=J\sum_{o\in \mathcal{O}(L)}\sum_{r\in
\sum(L_o)}(qa^{-1})^{{\rm
rot}(L_{o,r})}[L_o,r]R_{L_{o,r}}(q-q^{-1},a),
\end{eqnarray*}
where $J={1\over qa^{-1}+q^{-1}a}$.
\end{theorem}

\begin{remark} We have two remarks on Theorem \ref{Ja}.
\begin{enumerate}
\item[(1)] There is a coefficient $J$ in Theorem \ref{Ja}, for we use the
normalized Homflypt and Kauffman polynomials in this paper. So
$$J={\delta\over \mu}={{a-a^{-1}\over
q-q^{-1}}\over{a^2q^{-1}-a^{-2}q\over q-q^{-1}}+1}={1\over
qa^{-1}+q^{-1}a}.$$

\item[(2)] Since, for $L_o$ which contains a top inward crossing and any
$r\in \sum(L_o)$, we have $[L_o,r]=0$, there are some terms in the
right hand of Theorem \ref{Ja} which is equal to 0. In Theorem
\ref{Ja}, actually we can only consider balanced edge orientations
of $L$ which do not contain a top inward crossing. We add some 0
terms in the summation, which, you will see, is important to prove
our main Theorem \ref{rel}.
\end{enumerate}
\end{remark}

\subsection{Wu's formula}

In \cite{Wu}, Wu built a relation between the 3-variable KV
polynomial and the 3-variable MOY polynomial of 4-valent rigid
vertex spatial graphs. We shall only restrict ourselves to 4-valent
plane graphs.

\begin{figure}[htbp]
\centering
\includegraphics[width=5cm]{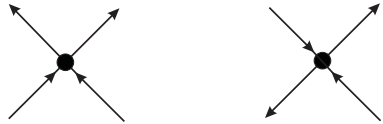}
\renewcommand{\figurename}{Fig.}
\caption{Balanced edge orientation near a vertex: the first is
 a crossing-like oriented vertex; the second is an alternatingly oriented vertex.}\label{cv}
\end{figure}

Let $G$ be a 4-valent plane graph. An edge orientation of $G$ is
\emph{balanced} if, at each vertex, among four edges incident with
the vertex, two edges are ``in" and two edges are ``out". Up to
rotation, there are two possible balanced edge orientations near a
vertex: one is crossing-like oriented and the other is alternatingly
oriented. (See Fig. \ref{cv}.)
Denote by $\mathcal{O}(G)$ the set of all balanced edge orientation
of $G$.

\begin{figure}[htbp]
\centering
\includegraphics[width=5cm]{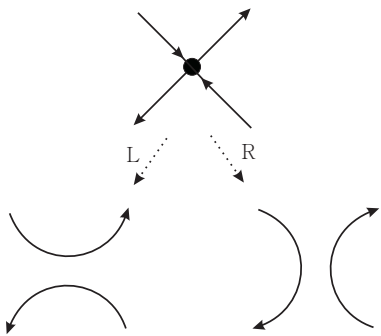}
\renewcommand{\figurename}{Fig.}
\caption{Two ways of smoothing an alternatingly oriented
vertex.}\label{spice-v}
\end{figure}

Given an $o\in \mathcal{O}(G)$, equipping $G$ with $o$, we obtain
$G_o$. To obtain oriented 4-plane graphs from $G_o$ we need to
splice all alternatingly oriented vertices. There are two ways of
smoothing an alternatingly oriented vertex of $G_o$ as shown in Fig.
\ref{spice-v}. A \emph{resolution} $r$ of $G_o$ is a choice of $L$
or $R$ smoothing of every alternatingly oriented vertex of $G_o$.
Denote by $\sum(G_o)$ the set of all resolutions of $G_o$. Given a
$r\in \sum(G_o)$, for each alternatingly oriented vertex $v$ of
$G_o$, we define the weight to be
\begin{displaymath}
[G_o,r;v] = \left\{
\begin{array}{ll}
-q & \textrm{if $r$ applies $L$ smoothing to the vertex $v$,}\\
-q^{-1} & \textrm{if $r$ applies $R$ smoothing to the vertex $v$.}
\end{array} \right.
\end{displaymath}
We can take the weight of crossing-like oriented vertex of $G_o$ to
be 1. The total weight $[G_o,r]$ of the resolution $r$ is defined to
be the product of weights of all vertices of $G_o$. Denote by
$G_{o,r}$ the oriented 4-valent plane graph obtained by applying $r$
to $G_o$.

\begin{theorem}(\textbf{Wu's Formula})  Let $G$ be an unoriented 4-valent plane graph. Then
\begin{eqnarray*}
[G]_D(q,q^{-1},a^2q^{-1})=J\sum_{o\in \mathcal{O}(G)}\sum_{r\in
\sum(G_o)}(qa^{-1})^{{\rm rot}(
G_{o,r})}[G_o,r][{G_{o,r}}]_R(q,q^{-1},a),
\end{eqnarray*}
where $J={1\over qa^{-1}+q^{-1}a}$.
\end{theorem}

\begin{remark}
Note that $[G_o,r;v]=-[D_e,z;v]$ in \cite{Wu}, since by comparing
Figs. \ref{He}, \ref{Ke} with Eqs. (1.11) and (1.9) in \cite{Wu}, we
know that the $[G]_R$ (resp. $[G]_D$) is equal to the result of
multiplying $(-1)^{|V(G)|}$ and $R(G)$ (resp. $P(G)$) in \cite{Wu}.
\end{remark}

\section{Two models and their relationship}

Now we are in a position to derive two oriented 4-valent graphical
models for the Kauffman polynomial. We call the $HJ$ model the
summation obtained by applying Jaeger's formula ($J$) firstly and
then Kauffman-Vogel model for the Homflypt polynomial ($H$).
Similarly, we call $WF$ model the summation obtained by applying
Kauffman-Vogel model for the Kauffman polynomial ($F$) and Wu's
formula ($W$). Let $L$ be an unoriented link diagram.

\begin{enumerate}
\item[(1)] The $HJ$ model: $D_L(q-q^{-1},a^2q^{-1})$\\
\begin{eqnarray*}
&=&J\sum_{o\in \mathcal{O}(L)}\sum_{r\in \sum(L_o)}(qa^{-1})^{{\rm
rot}(L_{o,r})}[L_o,r]\sum_{\
_{o,r}\!I}q^{i(_{o,r}I)-j(_{o,r}I)}[_{o,r}I]_R(q,q^{-1},a)\\
&=&J\sum_{o\in \mathcal{O}(L)}\sum_{r\in \sum(L_o)}\sum_{\
_{o,r}\!I}(qa^{-1})^{{\rm
rot}(_{o,r}I)}[L_o,r]q^{i(_{o,r}I)-j(_{o,r}I)}[_{o,r}I]_R(q,q^{-1},a),
\end{eqnarray*}
where the third summation runs over all oriented 4-valent plane
graphs: $_{o,r}I$'s, obtained from $L_{o,r}$ by applying $H$. The
second ``$=$"  holds since rot($L_{o,r}$)=rot($_{o,r}I$) for any
$_{o,r}I$.

\item[(2)] The $WF$ model: $D_L(q-q^{-1},a^2q^{-1})$\\
\begin{eqnarray*}
&=&\sum_{G}q^{i(G)-j(G)}J\sum_{o\in \mathcal{O}(G)}\sum_{r\in
\sum(G_o)}(qa^{-1})^{{\rm rot}(
G_{o,r})}[G_o,r][{G_{o,r}}]_R(q,q^{-1},a)\\
&=&J\sum_{G}\sum_{o\in \mathcal{O}(G)}\sum_{r\in
\sum(G_o)}q^{i(G)-j(G)}(qa^{-1})^{{\rm rot}(
G_{o,r})}[G_o,r][{G_{o,r}}]_R(q,q^{-1},a),
\end{eqnarray*}
where the first summation runs over all unoriented 4-plane graphs:
$G$'s, obtained from $L$ by applying $F$.
\end{enumerate}

Note that in the $HJ$ model, for different orientation $o$,
resolution $r$ and different ways of applying $H$ to $L_{o,r}$, the
obtained oriented 4-valent plane graphs: $_{o,r}I$'s, are all
different in the sense that the crossings and edges of $L$ are
labeled differently and kept unchanged after splicing some
crossings. Let $S$ be the set of all oriented 4-valent plane graphs
constructed from $L$ in the $HJ$ model. In other words, the terms in
the $HJ$ model are all different. However, in the $WF$ model, there
exist many terms whose corresponding $G_{o,r}$'s are the same. In
other words, for different $G$, different orientation $o$ and
resolution $r$, the obtained $G_{o,r}$'s are not all different.

\begin{figure}[htbp]
\centering
\includegraphics[width=8cm]{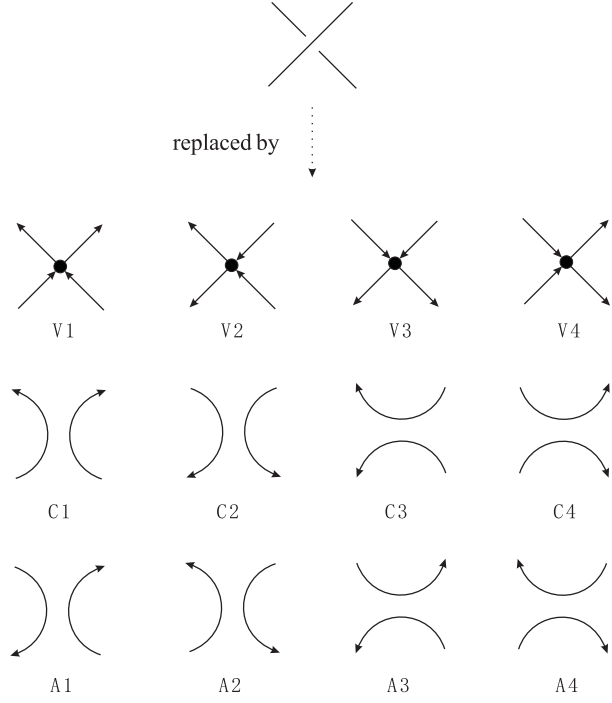}
\renewcommand{\figurename}{Fig.}
\caption{12 configurations used to replace an unoriented
crossing.}\label{rep-c}
\end{figure}

Furthermore, note that the set of different oriented 4-valent plane
in both models are the same: they are both the set of 4-valent plane
graphs obtained from $L$ by replacing each an unoriented crossing
$c$ of $L$ by one of the following twelve types of configurations:
$V_1$, $V_2$, $V_3$, $V_4$, $C_1$, $C_2$, $C_3$, $C_4$, $A_1$,
 $A_2$, $A_3$ and $A_4$ as shown in Fig. \ref{rep-c}. Of course, we demand
orientations of all local replacements are compatible on each edge
of the underlying 4-valent plane graph of $L$ when we construct
oriented 4-valent plane graphs from $L$.

Let $\mathcal{T}=\{G;o,r\}$ be the set of different ways: $G;o,r$'s,
of constructing oriented 4-valent plane graphs from $L$, which
correspond all terms in the $WF$ model. Now we rewrite the $HJ$ and
$WF$ models as the sum
\begin{eqnarray*}
\sum_{s\in S}c_s[s]_R {\ \ \  \textrm{and}\ \ \ } \sum_{G;o,r\in
\mathcal{T}}d_{G;o,r}[G_{o,r}]_R.
\end{eqnarray*}
Then we have

\begin{theorem}\label{rel} For each $s\in S$, there exists a subset
$T_s=\{G;o,r\in \mathcal{T}|s=G_{o,r}\}$ such that $T_{s_1}\cap
T_{s_2}=\emptyset$ for any $s_1,s_2\in s$, $\cup_{s\in
S}T_s=\mathcal{T}$ and $c_s=\sum_{G;o,r\in T_s}d_{G;o,r}$.
\end{theorem}

\noindent{\bf Proof.} We have shown the existence of $T_s$ such that
$T_{s_1}\cap T_{s_2}=\emptyset$ for any $s_1,s_2\in s$ and
$\cup_{s\in S}T_s=\mathcal{T}$ in the preceding several paragraphs.
Hence, it suffices for us to prove that $c_s=\sum_{G;o,r\in
T_s}d_{G;o,r}$. Let $s= _{o,r}I$. Recall that
\begin{eqnarray*}
c_s&=&J(qa^{-1})^{{\rm rot}(s)}[L_o,r]q^{i(_{o,r}I)-j({o,r}I)},\\
d_{G;o,r}&=&J(qa^{-1})^{{\rm rot}(G_{o,r})}q^{i(G)-j(G)}[G_o,r]\\
&=&J(qa^{-1})^{{\rm rot}(s)}q^{i(G)-j(G)}[G_o,r].
\end{eqnarray*}
Now let
\begin{eqnarray*}
c_{HJ}&=&[L_o,r]q^{i(_{o,r}I)-j(_{o,r}I)},\\
c_{WF}&=&\sum_{G;o,r\in T_s}q^{i(G)-j(G)}[G_o,r].
\end{eqnarray*}
Then we only need to prove $c_{HJ}=c_{FW}$. We suppose that $s$ is
the 4-valent plane graph obtained from $L$ by $v_i$ replacements of
the configuration $V_i$, $c_i$ replacements of the configuration
$C_i$ and $a_i$ replacements of the configuration $A_i$ for
$i=1,2,3,4$.

\begin{figure}[htbp]
\centering
\includegraphics[width=10cm]{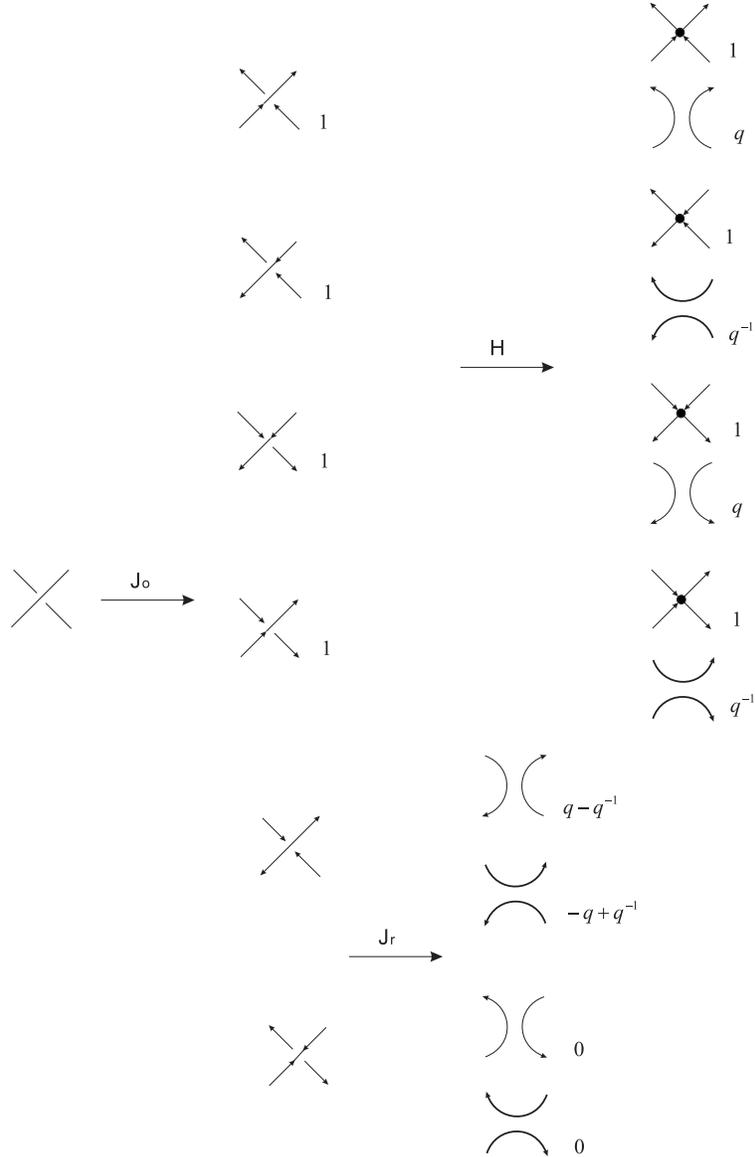}
\renewcommand{\figurename}{Fig.}
\caption{The $HJ$ model expansion, $J_o$ means orientation by
applying $J$, $J_r$ means resolution by applying $J$.}\label{HJ}
\end{figure}

By applying $J$ firstly then $H$, an unoriented crossing $c$ of $L$
will be replaced by one of eight oriented configurations firstly,
then each of the four crossing-like oriented configurations is
replaced by one of two configurations (see Fig. \ref{HJ}). The
corresponding weight in the product
$[L_o,r]q^{i(_{o,r}I)-j(_{o,r}I)}$ is labeled as the subscript in
that figure.

\begin{figure}[htbp]
\centering
\includegraphics[width=12cm]{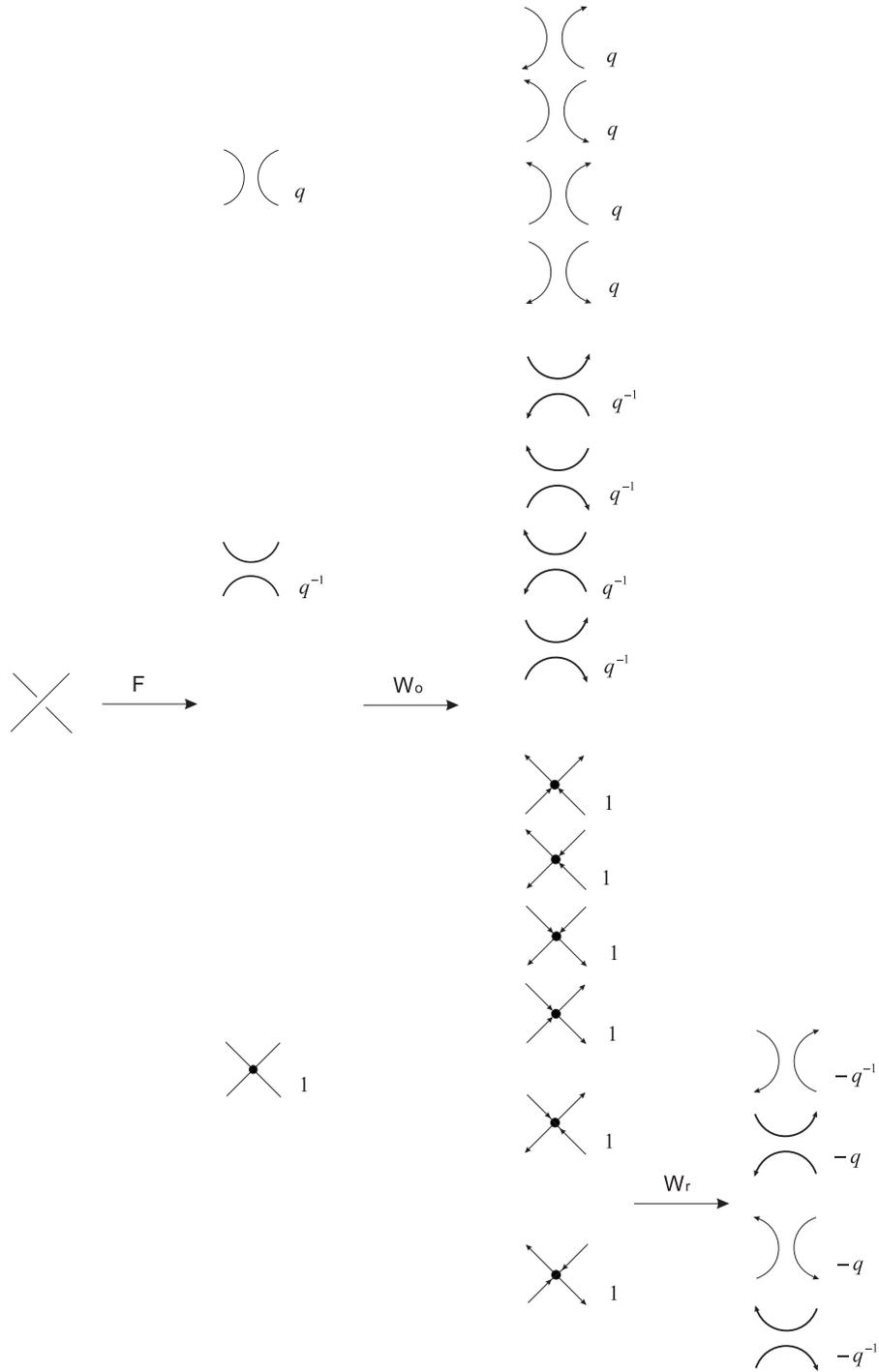}
\renewcommand{\figurename}{Fig.}
\caption{The $WF$ model expansion, $W_o$ means orientation by
applying $W$, $W_r$ means resolution by applying $W$.}\label{WF}
\end{figure}

Similarly, applying $F$ firstly then $W$, an unoriented crossing $c$
of $L$ will be replaced by one of three unoriented configurations:
$A$-smoothing, $B$-smoothing and the vertex replacement firstly,
then each of the two smoothings is replaced by one of four oriented
configurations and the vertex replacement is replaced by one of
eight oriented configurations (see Fig. \ref{WF}). The corresponding
weight in $q^{i(G)-j(G)}[G_o,r]$ is also labeled as the subscript in
the figure.

Note that in Figs. \ref{HJ} and \ref{WF}, the corresponding weights
of $V_1$, $V_2$, $V_3$, $V_4$, $C_1$, $C_2$, $C_3$ and $C_4$ in the
two figures are the same. Now we analyze the remaining four
configurations: $A_1$, $A_2$, $A_3$, $A_4$. There are two cases:

\noindent{\bf Case 1.} If $a_2\neq 0$ or $a_4\neq 0$, it is clear
that $c_{HJ}=0$. Now we consider $c_{WF}$. Without loss of
generality we suppose that $a_2\neq 0$. Note that $A_2$ appears
twice in Fig. \ref{WF}. This means we can obtain $a_2$ $A_2$
configurations by selecting $k$ $A_2$ configurations via applying
$F$ and then selecting the remaining $a_2-k$ $A_2$ configurations
via applying $W$ for any $k=0,1,\cdots,a_2$.  Note that
$\sum_{k=0}^{a_2}q^k(-q)^{a_2-k}=0$. Hence $c_{WF}=0$.

\noindent{\bf Case 2.} Otherwise, it means that $s$ does not contain
$A_2$ and $A_4$ configurations. Similarly, $A_1$ and $A_3$ appears
twice in Fig. \ref{WF}. Since
$(q-q^{-1})^{a_1}=\sum_{k=1}^{a_1}q^k(-q^{-1})^{a_1-k}$ and
$(-q+q^{-1})^{a_1}=\sum_{k=1}^{a_3}(-q)^{a_3-k}(q^{-1})^{k}$, we
have $c_{HJ}=c_{WF}$.

\noindent This completes the proof of Theorem \ref{rel}. $\Box$

\vskip0.3cm

Theorem \ref{rel} tells us that many terms of the $WF$ model add up
to one term of the $HJ$ model, hence, the $HJ$ model is more
efficient than $WF$ model. Now we provide an example to illustrate
Theorem \ref{rel}.

\begin{figure}[htbp]
\centering
\includegraphics[width=13cm]{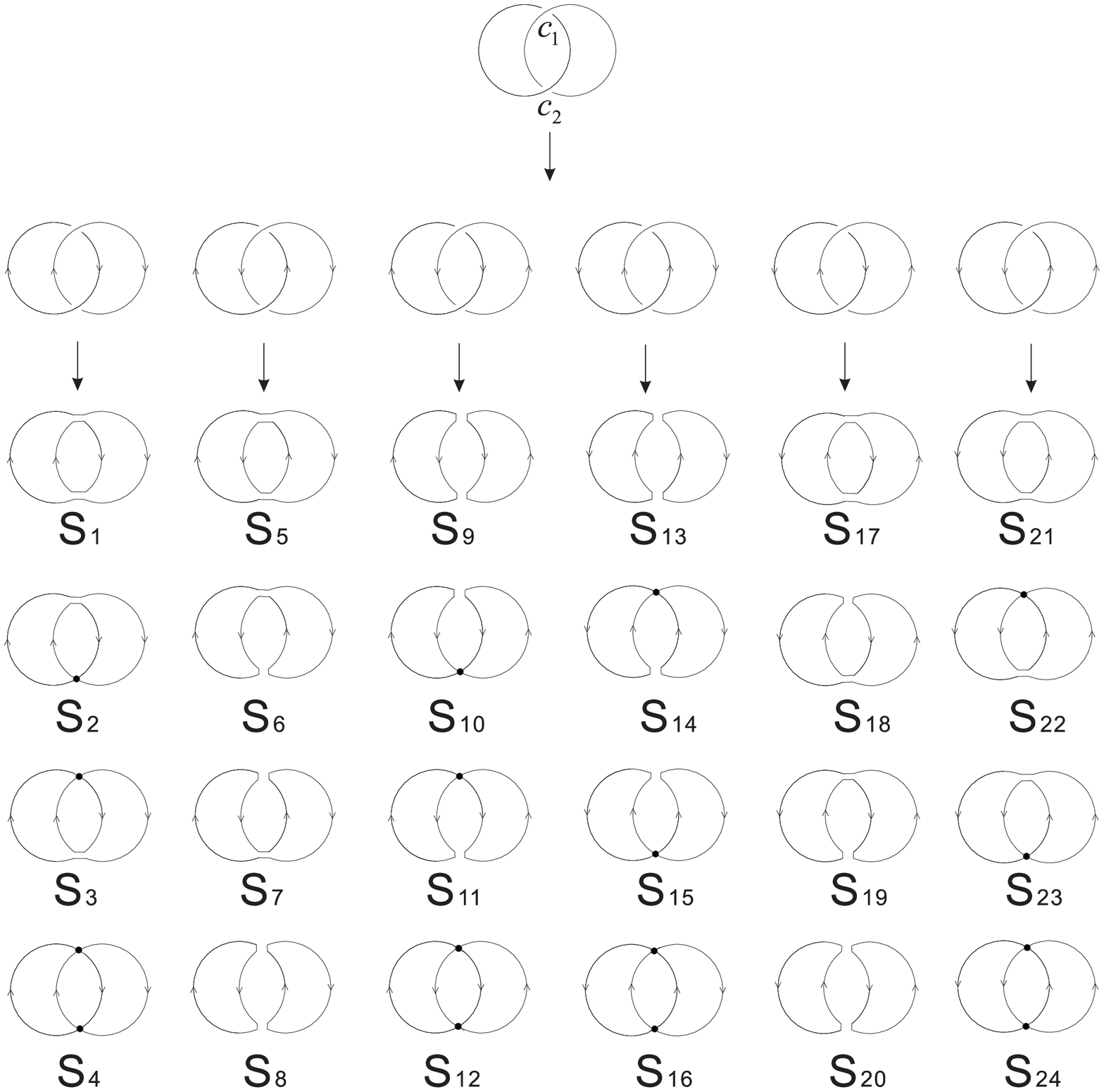}
\renewcommand{\figurename}{Fig.}
\caption{Hopf link, 6 balanced orientations and 24 states in the
$HJ$ model.}\label{hopf1}
\end{figure}

\begin{figure}[htbp]
\centering
\includegraphics[width=15cm]{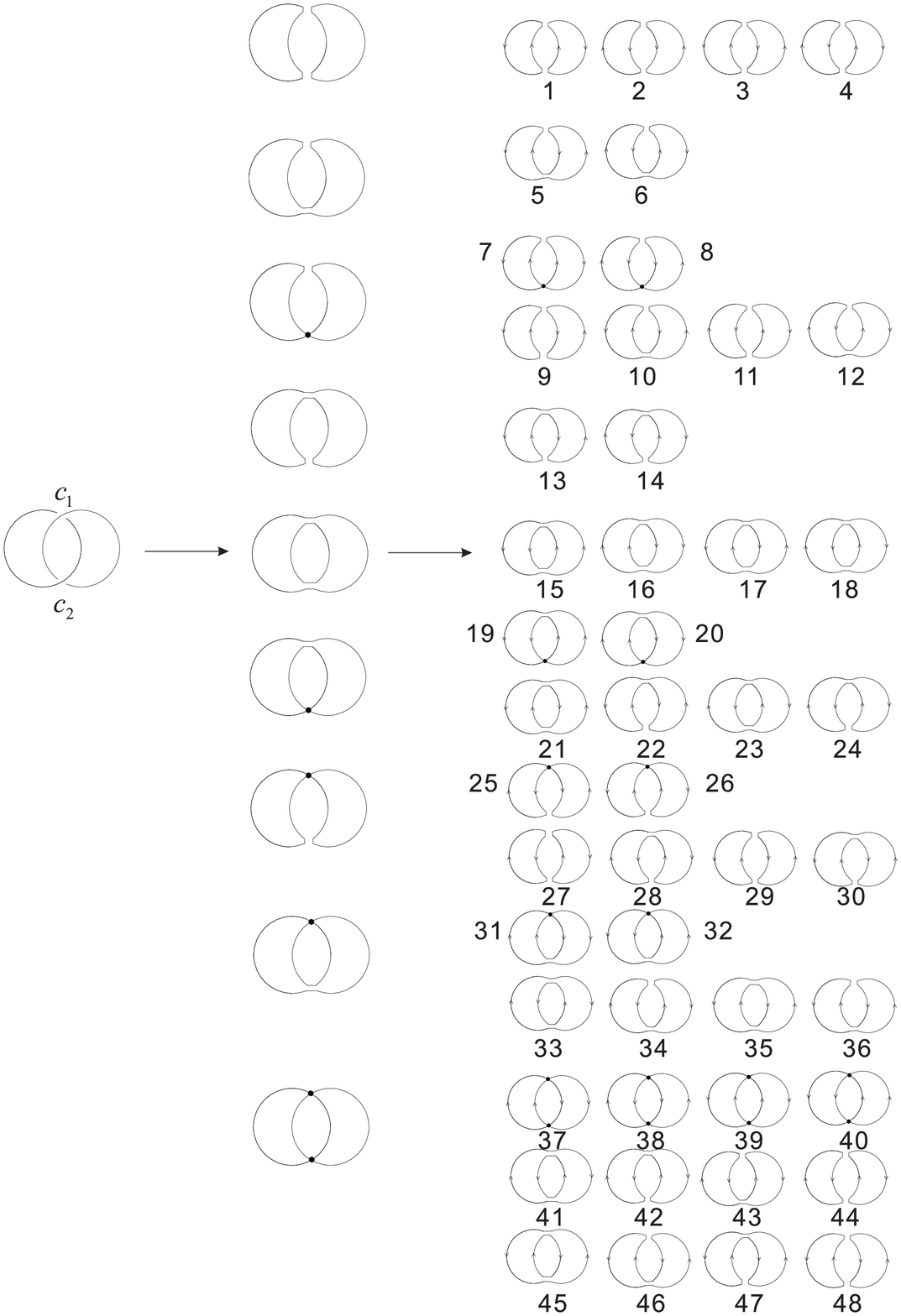}
\renewcommand{\figurename}{Fig.}
\caption{Hopf link, 9 corresponding unoriented 4-valent plane graphs
and 48 terms in the $WF$ model.}\label{hopf2}
\end{figure}

\begin{example}
The Hopf link.
\end{example}
We first expand the $R$ polynomial of the Hopf link based on $HJ$
model. There are six different balanced orientations for the Hopf
link, and twenty four oriented 4-valent plane graphs (i.e. states)
are constructed from the Hopf link (see Fig. \ref{hopf1}).
We then expand its $R$ polynomial based on $WF$ model. Nine
unoriented 4-valent plane graphs by applying $F$ are obtained
firstly, then forty eight terms all together are obtained by
applying $W$ to each 4-valent plane graph (see Fig. \ref{hopf2}).
For each $s_i$, its corresponding weights of the two crossings in
$HJ$ model, oriented 4-valent plane graphs corresponding to elements
of $T_{s_i}$ and the their corresponding weights of two crossings
are listed in Table 1. It is easy to verify that $c_{HJ}=c_{WF}$.

\begin{table}[htp]
\begin{center}
\begin{tabular}{|l|l|}\hline
  $s$ and its weight  & Graphs corresponding to elements of $T_s$ and their  \\
  $[w(c_1),w(c_2)]$   & corresponding weights of crossings $c_1$ and $c_2$  \\ \hline

  $s_{1}\ [q^{-1},q^{-1}]$  &  $16\ [q^{-1},q^{-1}]$                 \\ \hline
  $s_{2}\ [q^{-1},1]$       &  $20\ [q^{-1},1]$                      \\  \hline
  $s_{3}\ [1,q^{-1}]$       &  $31\ [1,q^{-1}]$                      \\ \hline
  $s_{4}\ [1,1]$            &  $37\ [1,1]$                           \\  \hline

  $s_5\ [q^{-1}-q,q^{-1}-q]$ &$18\ [q^{-1},q^{-1}]$, $23\ [q^{-1},-q]$, $33\ [-q,q^{-1}]$, $41\ [-q,-q]$  \\ \hline
  $s_6\ [q^{-1}-q,q-q^{-1}]$ &$14\ [q^{-1},q]$, $24\ [q^{-1},-q^{-1}]$, $28\ [-q,q]$, $42\ [-q,-q^{-1}]$  \\ \hline
  $s_7\ [q-q^{-1},q^{-1}-q]$ &$06\ [q,q^{-1}]$, $12\ [q,-q]$, $34\ [-q^{-1},q^{-1}]$, $43\ [-q^{-1},-q]$  \\ \hline
  $s_8\ [q-q^{-1},q-q^{-1}]$ &$04\ [q,q]$, $11\ [q,-q^{-1}]$, $27\ [-q^{-1},q]$, $44\ [-q^{-1},-q^{-1}]$  \\ \hline

  $s_9\ [q,q]$   &  $02\ [q,q]$  \\ \hline
  $s_{10}\ [q,1]$&  $08\ [q,1]$  \\ \hline
  $s_{11}\ [1,q]$&  $25\ [1,q]$  \\ \hline
  $s_{12}\ [1,1]$&  $38\ [1,1]$  \\ \hline

  $s_{13}\ [q,q]$&  $01\ [q,q]$\\ \hline
  $s_{14}\ [1,q]$&  $26\ [1,q]$ \\  \hline
  $s_{15}\ [q,1]$&  $07\ [q,1]$ \\ \hline
  $s_{16}\ [1,1]$&  $39\ [1,1]$  \\  \hline

  $s_{17}\ [0,0]$&  $17\ [q^{-1},q^{-1}]$, $21\ [q^{-1},-q^{-1}]$, $35\ [-q^{-1},q^{-1}]$, $45\ [-q^{-1},-q^{-1}]$\\ \hline
  $s_{18}\ [0,0]$&  $05\ [q,q^{-1}]$, $10\ [q,-q^{-1}]$, $36\ [-q,q^{-1}]$, $46\ [-q,-q^{-1}]$                    \\ \hline
  $s_{19}\ [0,0]$&  $13\ [q^{-1},q]$, $22\ [q^{-1},-q]$, $30\ [-q^{-1},q]$, $47\ [-q^{-1},-q]$                    \\ \hline
  $s_{20}\ [0,0]$&  $03\ [q,q]$, $09\ [q,-q]$, $29\ [-q,q]$, $48\ [-q,-q]$                                        \\  \hline

  $s_{21}\ [q^{-1},q^{-1}]$&  $15\ [q^{-1},q^{-1}]$  \\ \hline
  $s_{22}\ [1,q^{-1}]$&  $32\ [1,q^{-1}]$            \\ \hline
  $s_{23}\ [q^{-1},1]$&  $19\ [q^{-1},1]$            \\ \hline
  $s_{24}\ [1,1]$&  $40\ [1,1]$                      \\ \hline

  \end{tabular}
\vskip0.5cm {\footnotesize Table 1. $s$, $T_s$ and corresponding
weights.}
\end{center}
\end{table}

Now we simplify the $HJ$ model by deleting the 0-terms and obtain

\begin{theorem} Let $L$ be an unoriented link diagram. Then
\begin{eqnarray*}
D_L(q-q^{-1},a^2q^{-1})=J\sum_{\sigma}\left\{\prod_{c}w(\sigma,c)\right\}(qa^{-1})^{{\rm
rot}(\sigma)}[\sigma]_R(q,q^{-1},a),
\end{eqnarray*}
where $J={1\over qa^{-1}+q^{-1}a}$, the summation runs over all
oriented 4-valent plane graphs obtained from $L$ by replacing each
crossing by one of the $V_1$, $V_2$, $V_3$, $V_4$, $C_1$, $C_2$,
$C_3$, $C_4$, $A_1$, and $A_3$ (see Fig. \ref{rep-c}), the product
is over all crossings of $L$, and the weight $w(\sigma,c)$ depends
on the replacement of $c$ in $\sigma$ and is shown as the subscript
in Fig. \ref{HJ}.
\end{theorem}

\section{Trivalent graphical models}

There are also trivalent graphical models for both the $R$ and $D$
polynomial, see \cite{MOY,Fre} and \cite{CT} respectively. By
\emph{unoriented trivalent graph}, we mean a trivalent undirected
graph having two types of edges, ``thick" edges and ``common" edges
such that there is exactly one ``thick" edge incident to each
vertex. We take free loops to be special cases of unoriented
trivalent graphs. Clearly there is a many-to-one correspondence $f$
between such unoriented trivalent graphs and unoriented 4-valent
graphs and $f$ maps the unoriented trivalent graph to the unoriented
4-valent graph obtained from the unoriented trivalent graph by
contracting all thick edges. Under this correspondence and restrict
to planar graphs, by comparing the graphical calculus in this paper
with Eqs. (2.1)-(2.5) in \cite{CT}, you will see that the 3-variable
polynomial $P(G)$ in \cite{CT} and the 3-variable graph polynomial
$[f(G)]_D$ in \cite{KV} are completely the same. Clearly, for the
$D$ polynomial there is a bijection between terms of trivalent
graphical model \cite{CT} and terms of 4-valent graphical model
\cite{KV}. Hence the trivalent graphical model and the 4-valent
graphical model of the $D$ polynomial are essentially the same.

As for the $R$ polynomial, the relation between the trivalent
graphical model and the 4-valent graphical model is not very
immediate. As far as I know there is no trivalent graphical model
for the whole Homflypt polynomial, we only consider the special case
of the so-called Homflypt $n$-specializations, that is, we put
$z=q-q^{-1}$ and $a=q^n$ in the $R$ polynomial. In \cite{MOY},
Murakami, Ohtsuki and Yamada defined an invariant (we call it MOY
polynomial) of colored, oriented, trivalent plane graphs. In
\cite{Fre}, Freitas only considered a special case which only uses
colors 1 and 2 and called them \emph{classic} graphs. Note that
edges colored 1 correspond to ``common" edges and edges colored 2
correspond to ``thick" edges. In this special case, the
corresponding MOY polynomial is called the $\Gamma$-bracket.

\begin{theorem}\label{spe} Let $L$ be an oriented link diagram. Then
\begin{eqnarray}\label{Htr}
R_L(q-q^{-1},q^n)=q^{w(L)}{1\over
[n]}\sum_{G}(-q^{-1})^{s(G)}(-q)^{t(G)}<G>_n,
\end{eqnarray}
where $w(L)$ is the writhe of $L$,
$[n]=q^{-(n-1)}+q^{-(n-3)}+\cdots+q^{n-3}+q^{n-1}={q^n-q^{-n}\over
q-q^{-1}}$, the summation is over all $G$'s obtained from $L$ by
replacing each crossing by one of the two configurations shown as in
Fig. \ref{MOY-rep}, $s(G)$ (resp. $t(G)$) is the number of positive
(resp. negative) non-smoothed crossings of $L$ to form $G$, $<G>_n$
is the $\Gamma$-bracket of the classic graph $G$.
\begin{figure}[htbp]
\centering
\includegraphics[width=8cm]{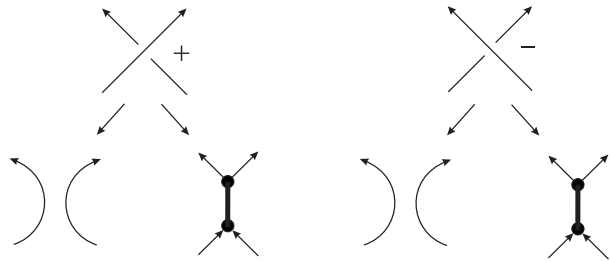}
\renewcommand{\figurename}{Fig.}
\caption{Two types of replacements to construct the classic
graphs.}\label{MOY-rep}
\end{figure}
\end{theorem}

\noindent{\bf Proof.} Theorem \ref{spe} is implicit in \cite{Fre}.
It can be obtained by combining Eqs. (2.12), (2.16) and (2.17) with
the fact that the $R$ polynomial in this paper is the normalized
regular invariant of the Homflypt polynomial. $\Box$

Now we simplify Eq. (\ref{Htr}) as follows.
\begin{eqnarray}
R_L(q-q^{-1},q^n)&=&q^{w(L)}\sum_{G}(-q^{-1})^{s(G)}(-q)^{t(G)}{1\over
[n]}<G>_n\nonumber\\
&=&q^{w(L)}\sum_{G}q^{t(G)-s(G)}{(-1)^{|V(G)|}\over [n]}<G>_n\nonumber\\
&=&\sum_{G}q^{i(G)-j(G)}{(-1)^{|V(G)|}\over [n]}<G>_n,\label{Htr1}
\end{eqnarray}
where $V(G)$ is the vertex set of $G$ and $i(G)$ (resp. $j(G)$) is
the number of positive (resp. negative) crossings smoothed to obtain
$G$.

Note that by contracting all ``thick" edges of a classic graph, we
obtain an oriented 4-valent plane graph. By putting $A=q$,
$B=q^{-1}$, $a=q^n$ in Theorem \ref{yy}, it is very similar to Eq.
(\ref{Htr1}). By comparing Eqs. (3.1)-(3.5) in \cite{Fre} and
graphical calculus for $[]_R$, it is not difficult for us to verify
that ${(-1)^{|V(G)|}\over [n]}<G>_n=[f(G)]_R(q,q^{-1},q^n)$ . We
leave the details to the readers. Note that for the Homflypt
$n$-specializations, the identity $(iii)'$ in Fig. \ref{id3} is
redundant.

Therefore, for the Homflypt $n$-specializations, the trivalent
graphical model in \cite{Fre} and the the 4-valent trivalent
graphical model are consistent. Clearly, for the $R$ polynomial
there is also a bijection between terms of trivalent graphical model
\cite{MOY,Fre} and terms of 4-valent graphical model \cite{KV}.

\vskip0.5cm

\noindent{\bf Acknowledgements}

\vskip0.2cm

This paper was completed during my visiting the Lafayette College. I
would like to thank Professor L. Traldi for introducing me to the
study of relations among various models of link polynomials and some
helpful conversations and comments. This work was also partially
supported by Grants from the National Natural Science Foundation of
China (No. 10831001) and the Fundamental Research Funds for the
Central Universities (No. 2010121007).





\bibliographystyle{model1b-num-names}
\bibliography{<your-bib-database>}

\begin{thebibliography}{00}

\bibitem{Alexander} J. W. Alexander, Topological invariants of knots and links, Trans. Amer. Math. Soc. 30 (1928) 275-306.

\bibitem{Bondy} J. A. Bondy and U. S. R. Murty, Graph theory with applications, The Macmillan press ltd, 1976.

\bibitem{BLM} R. D. Brandt, W. B. R. Lickorish and K. C. Millett, A polynomial invariant for unoirened knots and links, Invent. Math. 74 (1986) 563-573.

\bibitem{CT} C. Caprau, J. Tipton, The Kauffman polynomial and trivalent graphs, arXiv:1107.1210v2 [math.GT] 11 Jul 2011.

\bibitem{Carp} R. P. Carpentier, From planar graphs to embedded graphs-a new approach to Kauffman and Vogel's polynomial, J. Knot
Theory Ramifications 9(8) (2000) 975-986.

\bibitem{Conway} J. H. Conway, An enumeration of knots and links, and some of their algebraic properties, Computational Problems in Abstract Algebra,
Pergamon Press, New York (1970) 329-358.

\bibitem{Fer} E. Ferrand, On Legendrian knots and polynomial invariants, Proc. Amer. Amer. Soc. 130(4) (2001) 1169-1176.

\bibitem{Fre} N. R. B. Freitas, A combinatorial approach to the Homfly $n$-specializations, 2008.

\bibitem{Hom1} P. Freyd, D. Yetter, J. Hoste, W. B. R. Lickorish, K. Millett, and A. Ocneanu, A new polynomial invariant of knots and
links, Bull. Amer. Math. Soc. (N.S.) 12(2) (1985) 239-246.

\bibitem{Ho} C. F. Ho, A new polynomial invariant for knots and links-preliminary report, Abstracts Amer. Math. Soc. 6 (1985) 300.

\bibitem{Hug} S. Huggett, On tangles and matroids, J. Knot Theory Ramifications 14(7) (2005) 919-929.

\bibitem{Jae} F. Jaeger, Tutte polynomials and link polynomials, Proc. Amer. Math. Soc. 103 (1988) 647-654.

\bibitem{Jones} V. F. R. Jones, A polynomial invariant for knots via Von Neumann algebras, Bull. Amer. Math. Soc. 12 (1985) 103-111.

\bibitem{KB} L. H. Kauffman, State models and the Jones polynomial, Topology 26 (1987) 395-407.

\bibitem{KN} L. H. Kauffman, On knots, Annals of Mathematics Studies, No. 115, Princeton University Press, Princeton, New Jersey, 1987.

\bibitem{Kau1} L. H. Kauffman, Invariants of graphs in three-space, Trans. Amer. Math. Soc. 311(2) (1989) 697-710.

\bibitem{Kau2-variable} L. H. Kauffman, An invariant of regular isotopy, Trans. Amer. Math. Soc. 318(2) (1990) 417-471.

\bibitem{K} L. H. Kauffman, Knots and Physics, World Scientifc, 1991.

\bibitem{KV} L. H. Kauffman, P. Vogel, Link polynomials and a graphical calculus, J. Knot Theory Ramifications 1(1) (1992) 59-104.

\bibitem{Lik} W. B. R. Lickorish, Some link-polynomial relations, Math. Proc. Phil. Soc. 105 (1989) 103-107.

\bibitem{MOY} H. Murakami, T. Ohtsuki, S. Yamada, Homfly polynomial via an invariant of clored plane graphs, Enseign. Math. 44 (1998) 325-360.

\bibitem{Hom2} J. H. Przytycki, P. Traczyk, Invariants of links of Conway type, Kobe J. Math. 4 (1987) 115-139.

\bibitem{MBT} M. B. Thistlethwaite, A spanning tree expansion of the Jones polynomial,
Topology 26 (1987) 297-309.

\bibitem{Tutte} W. T. Tutte, A contribution to the theory of chromatic polynomials, Canad. J. Math. 6 (1954) 80-91.

\bibitem{Wu} H. Wu, On the Kauffman-Vogel and the Murakami-Ohtsuki-Yamada graph polynomials, arXiv:1107.5333v1
[math.GT] 26 July 2011.









\end{thebibliography}



\end{document}